\documentclass[11pt]{article}

\def\version{5.3.2017}\def\users{final-layout}  %
\usepackage{a4}
\usepackage{times,cite}
\usepackage[]{graphicx}
\usepackage{amsmath,amssymb,amsthm,mathrsfs,epsf,a4,color,graphicx,eucal}
\usepackage{ifthen}
\usepackage[normalem]{ulem}
\usepackage{mathrsfs}   
\usepackage{psfrag}
\usepackage{epsfig}
\ifthenelse{\equal{\users}{final-layout}}{}
{
\usepackage{fancyhdr}
\pagestyle{fancy}
\headheight=28pt
\rhead{\color{green}
Variational methods for Fick/Darcy flow...\\
by T.Roub\'\i\v{c}ek}
\chead{}
\lhead{\color{green}Version\,\version, file:\,\jobname
\\compiled:
\number\day.\number\month.\number\year\ at
\the\hour:\ifnum\minute<10 0\fi\the\minute\ h }
\usepackage{color}
\usepackage[notcite,notref,color]{showkeys}
\definecolor{labelkey}{rgb}{1.,.2,0.}

}
\newcount\hour \newcount\minute
\hour=\time
\divide \hour by 60
\minute=\time
\loop \ifnum \minute > 59 \advance \minute by -60 \repeat

\newcounter{myfigure}

\ifthenelse{\equal{\users}{final-layout}}{

	\newcommand{\COMMENT}[1]{}
	\newcommand{\DELETE}[1]{}

        \newcommand{\REM}[1]{\marginpar{\bfseries\tiny{\color{blue}}}}
         }
{
 
 \newcommand{\COMMENT}[1]{{\color{red}\uuline{#1}\color{black}}}
 \newcommand{\DELETE}[1]{{\color{brown}\sout{#1}\color{black}}}

 \newcommand{\REM}[1]{\marginpar{\bfseries\tiny{\color{blue}#1}}}
}

\newcommand\R{\mathbb R}

\newcommand\pl{\partial}
\renewcommand\d{\mathrm d}

\newcommand\DT[1]{\mathchoice
                 {{\buildrel{\hspace*{.1em}\text{\LARGE.}}\over{#1}}}
                 {{\buildrel{\hspace*{.1em}\text{\Large.}}\over{#1}}}
                 {{\buildrel{\hspace*{.1em}\text{\large.}}\over{#1}}}
                 {{\buildrel{\hspace*{.1em}\text{\large.}}\over{#1}}}}
\newcommand\DDT[1]{\mathchoice
   {{\buildrel{\hspace*{.1em}\text{\Large.\hspace*{-.1em}.}}\over{#1}}}
   {{\buildrel{\hspace*{.1em}\text{\large.\hspace*{-.1em}.}}\over{#1}}}
   {{\buildrel{\hspace*{.1em}\text{\large.\hspace*{-.1em}.}}\over{#1}}}
   {{\buildrel{\hspace*{.1em}\text{\large.\hspace*{-.1em}.}}\over{#1}}}}

\newcommand{\lineunder}[2]{\LU{\begin{array}[t]{c}\underbrace{#1}\vspace*{.5em}\end{array}}{\mbox{\footnotesize\rm #2}}}
\newcommand{\linesunder}[3]{\LSU{\begin{array}[t]{c}\underbrace{#1}\vspace*{.5em}\end{array}}{\mbox{\footnotesize\rm #2}}{\mbox{\footnotesize\rm #3}}}
\newcommand{\LU}[2]{\begin{array}[t]{c}#1\vspace*{-1em}\\_{#2}\end{array}}
\newcommand{\LSU}[3]{\begin{array}[t]{c}#1\vspace*{-1em}\\_{#2}\vspace*{-.5em}\\_{#3}\end{array}}
\newcommand{\GDir}{\Gamma_{\!\Dir}}  

\newcommand{\GNeu}{\Gamma_{\!\Neu}}

\newcommand{\Dir}{{\scriptscriptstyle\mathrm{D}}}  
\newcommand{\Neu}{{\scriptscriptstyle\mathrm{N}}}
\newcommand{\Cdot}{\hspace{-.1em}\cdot\hspace{-.1em}}
\newcommand{\Colon}{\hspace{-.15em}:\hspace{-.15em}}
\newcommand{\In}{\!\in\!}

\def\eq{\eqref}
\def\eps{\varepsilon}

\def\bbC{\mathbb C}

\def\bbM{\mathbb M}

\def\bbI{\mathbb I}
\def\bbK{\mathbb K}
\def\vecd{{\hspace*{-.15em}\vec{\hspace*{.15em}d}}}
\def\vece{{\hspace*{-.1em}\vec{\hspace*{.1em}e}}}
\def\BP{\noindent{\it Proof.\ }}
\def\EP{$\hfill\Box$\medskip}

\newenvironment{my-picture}[3]{\refstepcounter{myfigure}\label{#3}\setlength{\unitlength}{\textwidth}\begin{picture}(#1,#2)}{\end{picture}}

\theoremstyle{plain}

\newtheorem{theorem}{Theorem}[section]

\newtheorem{proposition}[theorem]{Proposition}

\newtheorem{remark}[theorem]{Remark}

\begin{document}
\begin{sloppypar}

\begin{center}
{\LARGE\bf Variational methods for steady-state Darcy/Fick flow\\in swollen and
poroelastic solids}

\bigskip
{\sc Tom\'a\v s Roub\'\i\v cek}

\bigskip

{
Mathematical Institute, Charles University,
Sokolovsk\'a 83, CZ-186~75~Praha~8,  Czech Republic,\\ and 
\\
Institute of Thermomechanics, CAS,
Dolej\v skova 5, CZ-182~00~Praha~8, Czech Republic.
}

\end{center}

\bigskip

\baselineskip=10pt

{\small
{\bf Abstract:} Existence of steady states in elastic media at small strains with diffusion 
of a solvent or fluid due to Fick's or Darcy's laws is proved by combining 
usage of variational methods inspired from static situations with Schauder's 
fixed-point arguments. In the plain variant, the problem consists in the force 
equilibrium coupled with the continuity equation, and the underlying operator
is non-potential and non-pseudomonotone so that conventional methods are not
applicable. In advanced variants, electrically-charged multi-component flows 
through an electrically charged elastic solid are treated, employing critical 
points of the saddle-point type. Eventually, anisothermal variants involving 
heat-transfer equation are treated, too. 
}

\bigskip


\baselineskip=12pt

\def\KONST{\kappa} 

\def\BIOT{M} 

\def\tr{{\rm tr}}

\def\Rsym{{\R_{\rm sym}^{d\times d}}}

\setcounter{myfigure}{0}

\section{Introduction}\label{sect_Intro}

\begin{figure}\begin{center}
\psfrag{s}{\small $\sigma$}
\psfrag{Fick}{\footnotesize\bf\begin{minipage}[t]{15em}materials undergoing swelling\\[.3em]\hspace*{3em}(Fick flow)\end{minipage}}
\psfrag{Darcy}{\hspace*{1em}\footnotesize\bf\begin{minipage}[t]{10em}porous materials\\[.3em]\hspace*{1em}(Darcy flow)\end{minipage}}
\psfrag{hydrid}{\tiny\bf HYDRIDE}
\psfrag{metal}{\tiny\bf METAL}
\psfrag{H-atom}{\scriptsize\bf H-atoms}
\psfrag{metal atoms}{\scriptsize\bf\begin{minipage}[t]{8em}metal\\\hspace*{-.0em}atoms\end{minipage}} 
\psfrag{pores}{\scriptsize\bf\begin{minipage}[t]{8em}pores
\\\hspace*{-.0em}filled 
\\\hspace*{-.0em}by water\end{minipage}} 
\psfrag{solvent}{\scriptsize\bf\begin{minipage}[t]{5em}solvent\\\hspace*{-.0em}unpacking\\\hspace*{-.0em}macro-\\\hspace*{-.5em}molecules\end{minipage}} 
\hspace*{-.0em}\includegraphics[width=35em]{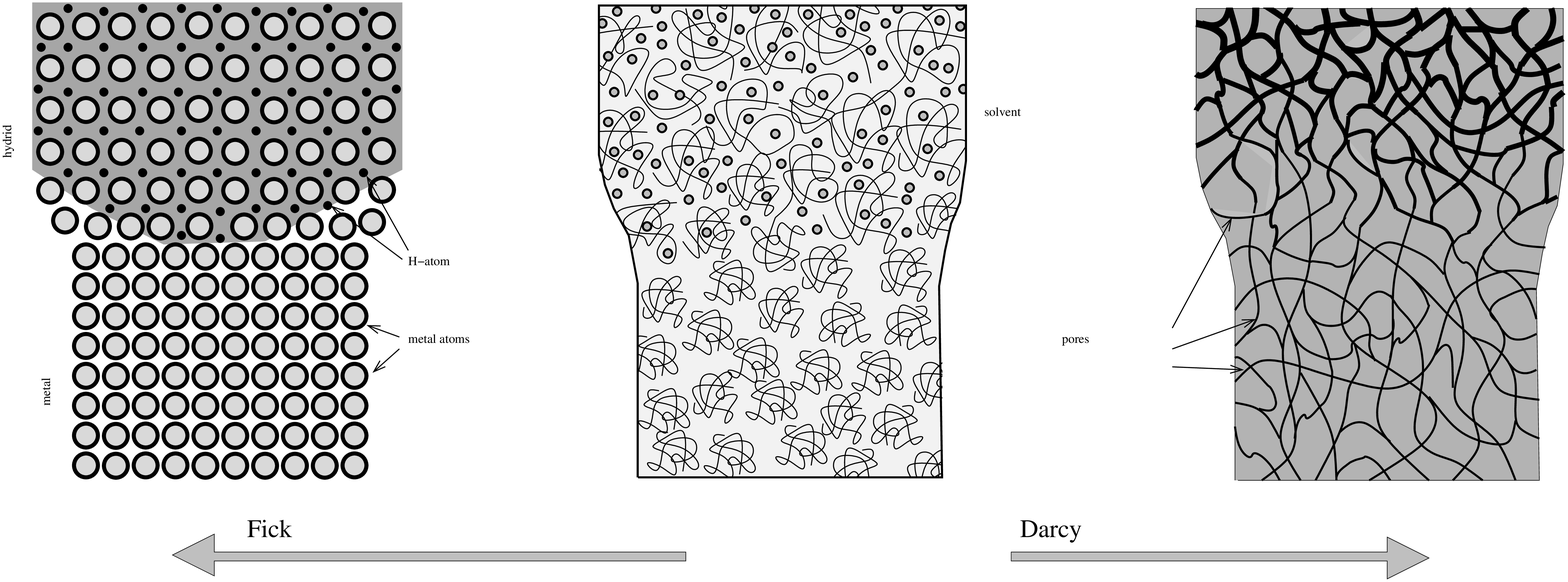}
\end{center}
\vspace*{-.6em}
\caption{Various underlying mechanisms coupling diffusion with volume 
expansion, ranging from an atomic via a molecular to a mesoscopical level:
\newline\hspace*{2em} Left: hydrogen diffusing inside a metallic atomic grid (metal/hydride transformation).
\newline\hspace*{2em} Middle: macromolecules in polymers unpacked by a diffusing solvent.
\newline\hspace*{2em} Right: water flowing through pores of a (poro)elastic rock or a concrete 
}\label{fig-variants-swelling}\end{figure}

Some elastic materials allow for a penetration of very small atoms into a 
solid atomic grid in crystalline metals or into spaces between big 
macromolecules of polymers. In the former case, the interstitial solute is 
hydrogen and such metals then undergo a so-called \emph{metal-hydride phase 
transformation} as schematically depicted in 
Fig.\,\ref{fig-variants-swelling}-left, see e.g.\ 
\cite{Latroche2004,Libowitz1994,RouTom14THSM}.
The latter mentioned mechanism occurs in polymers allowing for a diffusion of
a specific solvent causing unpacking of macromolecules as schematically 
depicted in Fig.\,\ref{fig-variants-swelling}-middle. In both cases, 
the solvent which diffuses thorough the 
elastic body influences the volume considerably (sometimes by tens of percents),
which is referred to as a {\it swelling}.
In turn, this swelling influences stress/strain distribution and therefore
also the diffusion process itself. The diffusion is driven rather
by the concentration gradient, which is referred to as a {\it Fick law}.

Other, microscopically different mechanism occurs in macroscopically solid 
materials that posses various pores or voids which are mutually connected and 
which allows for some fluids (sometimes referred to as a \emph{diffusant}) to 
flow thorough the solid. It is manifested macroscopically as a homogeneous 
mixture of a solid elastic body and fluid which diffuses throughout the 
volume as schematically depicted in Fig.\,\ref{fig-variants-swelling}-right. 
Examples are poroelastic rocks or porous polymers 
filled with water (and in the latter case possibly also with ionized hydrogen, 
i.e.\ protons, while the poroelastic polymer itself is negatively charged by
fixed dopands as used in polymer-electrolyte fuel cells \cite{ProWet09PEMF}). 
Interaction of 
solids with diffusants may be manifested by {\it squeezing}.
The diffusion is driven rather by the pressure gradient, which is referred to 
as {\it Darcy's law}.

The solid-diffusant interaction is surely a complicated multi-scale problem and 
a big amount of phenomenology is usually applied to build a simplified model.
A wide menagerie of models can thus be obtained, cf.\ \cite{Raja07HAMF} for a 
survey. Typically, small velocity of the diffusant is assumed. 
In this paper, we additionally assume small strains. 

In order to explain our ideas on a simple steady-state problem, we introduce a 
material model in which elasticity depends on an internal variable $c$ 
denoting a concentration of a diffusant or fluid. Concretely, in the simplest 
variant we will be concerned with the following boundary-value problem:
\begin{subequations}\label{EL-poro}\begin{align}
&{\rm div}\,\sigma+f=0&&\text{with }\ \ \ \sigma=\pl_e\varphi(e(u),c),&&
\\\label{EL-poro-mu}
&{\rm div}(\bbM(c)\nabla\mu)=0&&\text{with }\ \ \ \mu=\pl_c\varphi(e(u),c)&&
\end{align}
\end{subequations}
to be solved on a bounded Lipschitz domain $\Omega\subset\R^d$.
Later we will also investigate various generalization towards 
thermodynamics or coupling with electric field. The system \eq{EL-poro} should 
be completed by suitable boundary conditions, e.g.\ 
\begin{subequations}\label{poro-steady-BC}
\begin{align}\label{poro-steady-BC-u}
&u=0\ \text{ on }\GDir,\quad\ \sigma\vec{n}=g\ \,\text{ on }\,\GNeu,
\quad\ 
\\&\label{poro-steady-BC-c}
\bbM(c)\nabla\mu\cdot\vec{n}+\alpha\mu=\alpha\mu_{\rm ext}\ \text{ on }\Gamma,
\end{align}
\end{subequations}
where $\Gamma$ denotes the boundary of $\Omega$ divided into the parts 
$\GDir$ and $\GNeu$.

One should distinguish between the general \emph{steady-state}
(sometimes also called \emph{stationary})
situations and purely \emph{static}
situation (in the sense of steady states in {\it  thermodynamical equilibrium}). 
The former one means that 
all fields including the specific dissipation rates are independent of time, 
while the latter means in addition that the dissipation rates are zero (which
here means that all transportation processes vanish).

The plan of this paper is as follows: In Section~\ref{sec-FE},
we start with examples of free energies $\varphi$ typically arising in 
problems from continuum mechanics of poroelastic or swelling-exhibiting
materials. The static case enjoys a full variational structure, which is 
presented in Section~\ref{sec-static}. It can be exploited for general
steady-state problems where it is combined with the Schauder fixed point
technique, which is scrutinized in Section~\ref{sec-steady}. The steady-state 
problems \eq{EL-poro}-\eq{poro-steady-BC} suggest interesting thermodynamically
consistent anisothermal augmentation by a heat-transfer equation. This again 
combines variational technique for auxiliary static-like problems with 
Schauder's fixed point arguments, as presented in Section~\ref{sec-temp}. In
each Sections~\ref{sec-static}--\ref{sec-temp}, we also consider a (possibly
electrically charged) multi-component generalization which leads to a
necessity of considering critical points of a saddle-point-type variational
problems instead of mere minimizers. Various remarks close this paper in 
Section~\ref{sec-rem}, addressing in particular uniqueness issue, 
electroneutrality limit, or a general evolutionary thermodynamical context.

\section{Convex free energies leading to Fick or Darcy flows}\label{sec-FE}

With a certain simplification, one may say that two basic alternatives in the 
isotropic continuum are that either the strain is directly influenced by a 
concentration $c$ of the solvent which diffuses through the material
or the solvent has its own pressure which influences the overall stress and 
eventually indirectly the strain too. 
These two options microscopically reflect that the solvent does not 
directly contribute to the pressure but can change the reference 
configuration by influencing the atomic grid (as in the mentioned 
metal-hydride transformation) or the the solvent flows
though pores under its pressure and, from the macroscopical viewpoint, 
the resulted porous solid is phenomenologically homogenized.

The simplest form of the free energy in the former option is 
\begin{subequations}\label{swelling-strain}
\begin{align}\label{free-energy-swelling-strain}
&\hspace*{-3em}\!\!\!\!\varphi(e,c)=\frac12\bbC e^{}_{\rm el}\Colon e^{}_{\rm el}+
\KONST\Big(c\big({\rm ln}\frac c{c_{\rm eq}}{-}1\big)+\delta_{(0,\infty)}^{}(c)\Big)
\ \text{ with }\ e^{}_{\rm el}:=e{-}
Ec{+}Ec_{\rm eq}^{},
\end{align}
with $\bbC\in\R^{d\times d\times d\times d}$ a symmetric positive definite
4rt-order elasticity tensor, $\KONST\ge0$ a coefficient weighting chemical
versus mechanical effects,
$E\in\Rsym$ a matrix of the swelling coefficients
(which, in isotropic materials, is an identity matrix up to a coefficient)
and with $c_{\rm eq}$ an equilibrium concentration minimizing
$\varphi(0,\cdot)$. Moreover, $\delta_S(c)$
denotes the indicator function, i.e.\ $\delta_S(c)=0$ if $c\In S$ while
$\delta_S(c)=\infty$ otherwise. This yields the stress
\begin{align}\label{swelling-strain-s}
&\sigma=\pl_e\varphi(e,c)=\bbC e^{}_{\rm el}=\bbC(e{-}Ec{+}Ec_{\rm eq}^{})\ 
\intertext{and the chemical potential 
$\mu=\pl_c\varphi(e,c)=\KONST{\rm ln}\,c-E\Colon\sigma$ for $c>0$, or, taking 
into account that $\varphi(e,\cdot):\R\to\R\cup\{\infty\}$ is a proper convex 
function which is nonsmooth, rather}
\label{swelling-strain-mu}
&\mu\in\pl_c^{}\varphi(e,c)=
\begin{cases}
\displaystyle{\big\{\KONST{\rm ln}\,(c/c_{\rm eq}^{})
-E\Colon\sigma\big\}}&\text{if $c>0$},
\\\qquad\emptyset&\text{if $c\le0$};\end{cases}
\end{align}\end{subequations}
note that the last term $E\Colon\sigma$ is the pressure. Also note that 
$\lim_{c\to0+}\pl_c^{}\varphi(e,c)=-\infty$, which causes that indeed
$\pl_c^{}\varphi(e,0)=\emptyset$, cf.\ also \cite[Fig.\,11]{Roub13NPDE}.

The latter option relies on the idea that the diffusant (fluid) fully occupies 
the pores (i.e.\ so-called
saturated flow) whose volume is proportional to $\tr\,e(u)$ 
through a coefficient $\beta>0$ 
and its pressure is 
\begin{align}\label{swelling-p-fluid}
p_{\rm fld}^{}=\BIOT\big(\beta\tr\,e(u)-c+c_{\rm eq}^{}\big)
\end{align} 
with $\BIOT$ allowing an interpretation as the compressibility modulus of the 
fluid; actually, $M$ and $\beta$ are the so-called \emph{Biot modulus}
and the \emph{Biot coefficient} used in 
conventional models of porous media \cite{Biot41GTTS}. 
The (positive) parameter 
$c_{\rm eq}^{}$ denotes the equilibrium concentration and is considered here 
fixed, being related with the porosity of the material
which is here considered as a fixed material property.
This pressure is then summed up with the stress $\sigma_{\rm el}$ in the elastic 
solid. In isotropic materials, the total stress is then 
$\sigma=\sigma_{\rm el}+\beta p_{\rm fld}^{}\bbI
$. 
In such simplest variant, it leads to the potential 
\begin{subequations}\label{swelling-stress}
\begin{align}\label{free-energy-swelling-stress}
\hspace*{-4em}\varphi(e,c)=\frac12\bbC e\Colon e+
\frac12\BIOT(\beta\tr\,e{-}c{+}c_{\rm eq}^{})^2
+\KONST\Big(c\big({\rm ln}\frac c{c_{\rm eq}^{}}-1\big)
+\delta_{(0,\infty)}^{}(c)\Big),
\end{align}
which yields the stress and the chemical potential
\begin{align}\label{swelling-stress-s}
&\sigma=\pl_e\varphi(e,c)=\bbC e+
\beta\BIOT(\beta\tr\,e{-}c{+}c_{\rm eq}^{})\bbI\ \ \ \text{ and }\ \ \ 
\\\label{swelling-stress-mu}
&\mu\in\pl_c^{}\varphi(e,c)=
\begin{cases}\big\{\BIOT(c{-}\beta\tr\,e{-}c_{\rm eq}^{})
+\KONST{\rm ln}(c/c_{\rm eq}^{})\big\}&\text{if $c>0$},
\\\qquad\emptyset&\text{if $c\le0$},\end{cases}
\end{align}\end{subequations}
the emptyness of the subdifferential for $c\le0$
being analogous to \eqref{swelling-strain-mu}.

Then, choosing still a standard ansatz $\bbM(c)=c\bbM^{}_0$, the 
flux $j=-\bbM(c)\nabla\mu$
turns into 
\begin{align*}
\!\!j
=\!\!\!\!\!\linesunder{-\,c\bbM^{}_0\nabla p}{Darcy}{law}\!\!\!\!\!\!\!\!\!
\linesunder{-\,\bbM^{}_0\KONST\nabla c}{Fick}{law}
\ \text{ with }\ p=\begin{cases}
E\Colon\bbC(e{-}Ec{+}Ec_{\rm eq}^{})\!\!\!&\text{in case }\eq{swelling-strain},
\\[-.3em]
\BIOT(\beta\tr\,e{-}c{+}c_{\rm eq}^{})\!\!
&\text{in case }\eq{swelling-stress}
\end{cases}
\end{align*}
provided $c>0$. Depending on $\KONST$, either Darcy's mechanism or the 
Fick's one may dominate. Note also that $|p|=\mathscr{O}(|E|)$ in the case 
\eq{swelling-strain}.
An interesting phenomenon is that the equilibrium concentration 
$c_{\rm eq}^{}$ does not influence $\nabla\mu$ and can influence the
solution only through the boundary conditions \eq{poro-steady-BC}.
The mass-conservation equation \eq{poro-steady2}
reveals that the pressure gradient $\nabla p$ is needed and it also reveals
a certain ``optical'' difficulty because there is no obvious estimate on 
$\nabla p$. Indeed, in the evolution variant (except simple linearized 
problems like \eq{free-energy-swelling-stress} for $\KONST=0$ and 
$\bbM(\cdot)$ constant, like considered in fact in \cite{ShoSte04DPEM}),
a certain ``regularization'' of the problem 
seems to be necessary, by introducing a suitable phase field and its 
gradient as in \cite{RouTom14THSM}, or a gradient of $c$ leading to 
the Cahn-Hilliard ``capillarity-like'' model, or a gradient of $e(u)$
leading to a 2nd-grade nonsimple material concept. 

Yet, the standard definition of a weak solution to the boundary-value problem 
\eq{EL-poro}--\eq{poro-steady-BC} avoids explicit occurrence of $\nabla p$,
and it indeed works if $\bbM$ is constant; note that the fixed-point
argument used in the proofs of all ``non-static'' 
Propositions~\ref{prop-steady}--\ref{prop-steady-el-therm}
becomes rather trivial because the distribution of the 
chemical potential $\mu$ is then fully determined by $\mu_{\rm ext}$ in 
the boundary condition \eq{poro-steady-BC-c}. 
However, quite surprisingly, the non-static steady-state models allow for 
a lot of results without any regularization even when $\bbM$ depends on 
$c$; cf.\ also Remark~\ref{rem-other-FP} below. 

The following standard notation will be used: 
$L^p$ will stand for the Lebesgue
spaces of measurable functions whose $p$-power is integrable and $W^{1,p}$ 
for Sobolev spaces whose distributional derivatives 
are in $L^p$-spaces. For $p=2$, we abbreviate $H^1=W^{1,2}$.
Moreover, we use the standard notation  $p'=p/(p{-}1)$, and 
$p^*$ for the Sobolev exponent $p^*=pd/(d{-}p)$ for $p<d$ while
$p^*<\infty$ for $p=d$ and $p^*=\infty$ for $p>d$,
and the ``trace exponent'' $p^\sharp$ defined 
as $p^\sharp=(pd{-}p)/(d{-}p)$ for $p<d$ while
$p^\sharp<\infty$ for $p=d$ and $p^\sharp=\infty$ for $p>d$.
Thus, e.g., $W^{1,p}(\Omega)\subset L^{p^*}\!(\Omega)$ or 
$L^{{p^*}'}\!(\Omega)\subset W^{1,p}(\Omega)^*$. In the vectorial case,
we will write $L^p(\Omega;\R^n)\cong L^p(\Omega)^n$ 
and $W^{1,p}(\Omega;\R^n)\cong W^{1,p}(\Omega)^n$. 
We will also use the notation 
$H_\Dir^1(\Omega;\R^d):=\{u\In H^1(\Omega;\R^d);\ u|_{\GDir}=0\}$.
Eventually, ``$\,\cdot\,$''
or ``$\,:\,$'' denotes the scalar products of vectors or matrices, respectively.

Moreover, in what follows we will use the standard notation ``$\pl$'' either
for the partial derivative or a (partial) convex subdifferential.  
Without restricting generality towards applications of our
interest, we assume $\varphi$ nonsmooth in terms of $c$ only so that
$\pl_e^{}\varphi$ will be single-valued. Our variational methods 
will need coercivity of the involved convex or convex/concave functionals but 
not necessarily their controlled growth, so the related Euler-Lagrange
equations are to be understood in a variational rather than a conventional 
weak sense. Avoiding the growth conditions is one of the benefits of the
variational approach to this problem.

\section{Problems in thermodynamical equilibrium}\label{sec-static}

Let us now focus on special steady-state problems where the 
 dissipation rate is not only constant in time but just zero. Then also 
temperature is not influenced by the mechanical/diffusion part and
we can ignore the heat transfer. As we already said, such problems are called 
static. The diffusion is related with the dissipation rate (and entropy and 
heat production rate), namely that the overall dissipation (or also the 
heat-production) rate is $\int_\Omega\bbM\nabla\mu\Cdot\nabla\mu\,\d x$, 
cf.\ also to \eq{poro-therm-engr-} below. 
This implies $\nabla\mu$ everywhere on $\Omega$, and in particular 
also the boundary flux $j=-\bbM\nabla\mu$ is to be zero.
When assuming $\Omega$ connected in such static case, 
$\nabla\mu=0$ on $\Omega$ leads to that $\mu$ is constant. 
Let us denote this constant by $\bar\mu$. 

Therefore, solvability of such problem essentially requires either the system 
to be in equilibrium with the external environment or to be isolated. The 
former option is rather trivial: assuming $\mu_{\rm ext}=\bar\mu$ in 
\eq{poro-steady-BC} with a given constant $\bar\mu$ and $\alpha>0$, it fixes 
$\mu=\bar\mu$ and then one can eliminate $c$. From $\bar\mu=\pl_c\varphi(e,c)$, 
we can then find the concentration 
$c=[\pl_c\varphi(e,\cdot)]^{-1}(\bar\mu)$ as a function of $e$.
If $\varphi(e,\cdot)$ is convex,
we can even write a bit more specifically 
\begin{align}\label{swelling-small-c=c(q)}
c=\pl_\mu\varphi^*(e,\bar\mu)
\end{align}
with $\varphi^*(e,\cdot)$ denoting the convex conjugate function of 
$\varphi(e,\cdot)$. Note that, even if $\varphi(e,\cdot)$ is not smooth as in 
the examples in Section~\ref{sec-FE}, $\varphi^*(e,\cdot)$ is indeed 
single-valued if the natural requirement $\pl_{cc}^2\varphi>0$ holds.
The concentration $c$ can thus be completely eliminated. 

The latter option (i.e.\ the boundary permeability coefficient $\alpha=0$) is 
more interesting. When the profile $e=e(u(x))$ is known, the overall amount of 
solvent $\int_\Omega c\,\d x=C^{}_{\rm total}$ depends monotonically on $\bar\mu$ 
due to \eq{swelling-small-c=c(q)}, which allows us to specify $\bar\mu$ if 
$C^{}_{\rm total}$ is given. Yet, the displacement $u$ is a part of solution 
itself so that, unfortunately, it does not seem possible to fix $\bar\mu$ just
from knowing $C^{}_{\rm total}$.
In this isolated situation, it is natural to prescribe the total amount of 
diffusant
\begin{align}\label{swelling-small-constraint}
\int_\Omega c\,\d x=C^{}_{\rm total}\qquad\text{with\ \ \ }C^{}_{\rm total}\ge0
\text{\ \ \ given}.
\end{align}
Interestingly, this special steady-state (=static) case enjoys 
a full variational structure at least in the sense that some (if not all)
solutions can be obtained by such way as critical points. 
%
The general steady-state system \eq{EL-poro} then modifies to 
\begin{subequations}\label{EL-poro-stat}\begin{align}\label{poro-steady1}
&{\rm div}\,
\pl_e\varphi(e(u),c)+f=0,
\\\label{poro-steady2}
&\pl_c\varphi(e(u),c)\ni\bar\mu=\text{ some constant},
\end{align}
\end{subequations}
to be coupled with \eq{swelling-small-constraint} and with the boundary 
conditions \eq{poro-steady-BC-u}. Let us recall that 
$\varphi(e,\cdot)$ is allowed to be nonsmooth so that $\pl_c\varphi$
maybe set-valued so that \eq{poro-steady2} is an inclusion rather than
equation.
More specifically, 
the mentioned variational structure consists in the following constrained 
minimization problem is of a certain relevance:
\begin{align}\label{swelling-small-static-funct-u-c}
\left.\begin{array}{ll}
\text{Minimize }&\ \displaystyle{(u,c)\mapsto
\int_\Omega\varphi(e(u),c)-f\Cdot u\,\d x-\int_{\GNeu}\!\!g\Cdot u\,\d S} 
\\\text{subject to }&\displaystyle{\int_\Omega c\,\d x=C^{}_{\rm total}},
\ \ \ \ u\In H_\Dir^1(\Omega;\R^d),\ \ \ \ c\In L^1(\Omega).\ \ 
\end{array}\right\}
\end{align}

\begin{proposition}[Existence 
of static solutions]\label{prop-static}
Let $\varphi:\R_{\rm sym}^{d\times d}\times\R\to\R\cup\{\infty\}$ be convex,
lower semicontinuous, and coercive in the sense that 
$\varphi(e,c)\ge\epsilon|e|^2+\epsilon|c|^{1+\epsilon}$ for some 
$\epsilon>0$, 
${\rm meas}_{d-1}(\GDir)>0$,
$f\In L^{{2^*}'}\!(\Omega;\R^d)$, and $g\In L^{{2^\sharp}'}\!\!(\GNeu;\R^d)$.
Moreover, let the body be isolated (i.e.\ $\alpha=0$) and 
the overall content $C^{}_{\rm total}$ be given, assuming 
$\varphi(0,C^{}_{\rm total}/{\rm meas}_d(\Omega))<\infty$.
Then
the boundary-value problem \eq{EL-poro-stat} with 
boundary conditions \eq{poro-steady-BC-u}
and with the side condition \eq{swelling-small-constraint}
possesses 
at least one variational solution 
$(u,c)\in H_\Dir^1(\Omega;\R^d)\times L^{1+\epsilon}(\Omega)$ whose chemical 
potential 
is constant over $\Omega$ in the sense that 
there exists a constant $\bar\mu$ such that $\int_\Omega\varphi(e(u),c)
-\bar\mu c\,\d x\le\int_\Omega\varphi(e(u),\tilde c)
-\bar\mu\tilde c\,\d x$ for all $\tilde c\in C(\Omega)$, in particular 
\eq{poro-steady2} holds 
a.e.\ on $\Omega$ if $\pl_c^{}\varphi(e(u),c)$ exists.
\end{proposition}

\BP
We use the direct method for minimization of the convex coercive functional 
subject to the affine constraint \eq{swelling-small-static-funct-u-c}.
Note that the assumption $\varphi(0,C^{}_{\rm total}/{\rm meas}_d(\Omega))<\infty$
guarantees that this problem is feasible, i.e.\ its admissible set
contains at least one pair $(u,c):=(0,C^{}_{\rm total}/{\rm meas}_d(\Omega))$.
We thus obtain a minimizer 
$(u,c)\in H_\Dir^1(\Omega;\R^d)\times L^{1+\epsilon}(\Omega)$.

As the functional in \eq{swelling-small-static-funct-u-c} is convex, lower
semicontinuous, and 
the constraint is affine, introducing the Lagrange multiplier $\bar\mu$ to the 
scalar-valued constraint \eq{swelling-small-constraint} which is involved in 
\eq{swelling-small-static-funct-u-c}, this problem is equivalent to finding 
a critical point of the  functional (a so-called Lagrangian):
\begin{align}\label{swelling-small-static-funct-u-c-mu}
\hspace*{-4em}!\!\!(u,c,\bar\mu)\mapsto J(u,c,\bar\mu):=
\int_\Omega\varphi(e(u),c)-f\Cdot u-\bar\mu c
\,\d x-\!\int_{\GNeu}\!g\Cdot u\,\d S+C^{}_{\rm total}\bar\mu.\!\!
\end{align}
If this functional is smooth, from putting the G\^ateaux derivative 
with respect to $c$ zero, we can read that $\pl_c\varphi(e(u),c)=\bar\mu$ 
on $\Omega$, i.e.\ the chemical potential is constant. In a general nonsmooth 
case, by disintegration of the condition 
$\pl_c^{}J(u,c,\bar\mu)\ni0$ we obtain the inclusion \eq{poro-steady2} a.e.\ on 
$\Omega$. From putting $\pl_{\bar\mu}^{}J(u,c,\bar\mu)=0$,
we can read that 
$\int_\Omega -c\,\d x+C^{}_{\rm total}=0$, i.e.\
that the affine constraint \eq{swelling-small-constraint} 
is satisfied. Eventually, if the G\^ateaux derivative 
$\pl_u^{}J(u,c,\bar\mu)$ 
exists, it must be zero, and 
we can read the equilibrium equation \eq{poro-steady1} with the 
boundary conditions \eq{poro-steady-BC-u}
in the weak formulation. In a general case,  
\eq{poro-steady1}--\eq{poro-steady-BC-u} holds formally in a variational 
sense, i.e.\ 
the solution $u$ is considered as a minimizer of $J(\cdot,c,\bar\mu)$.
\EP

The above proof reveals the role of the scalar $\bar\mu$ as the Lagrange 
multiplier to the scalar constraint $\int_\Omega c\,\d x=C^{}_{\rm total}$,
and it is a vital part of the solution. This is consistent with the 
above observation that, if the mechanical part of the solution
$u$ is fixed, then $\bar\mu$ is determined by $C^{}_{\rm total}$.

If $\varphi$ is  strictly convex (as e.g.\ in the examples from 
Sect.\,\ref{sec-FE}) the minimization problem 
\eq{swelling-small-static-funct-u-c} has a unique solution,
although the relation to solutions of the static problem in question
may be more delicate, cf.~Remark~\ref{rem-uniqueness} below. 

An interesting and useful generalization of the above basic scenario
is towards a {\it multi-component fluid} with $N\ge2$ components which 
can even be {\it electrically charged} with specific 
charges $z=(z_1,...,z_N)$. Also the elastic medium can charged by some
dopands with the specific charge $z_{_{\rm DOP}}$. In the static problems 
there are no chemical reactions. The scalar-valued chemical potential is now 
to be replaced by an $\R^N$-valued electro-chemical potential 
$$
\mu=\pl_c^{}\varphi(e,c)+z\phi
$$ 
with $\phi$ the electrostatic potential. In static problems, again $\mu$ 
constant. The system \eq{EL-poro-stat} then augments to
\begin{subequations}\label{poro-stat-el}
\begin{align}\label{poro-stat-el1}
&{\rm div}\,\pl_e\varphi(e(u),c)+f=z_{_{\rm DOP}}\nabla\phi
&&\text{on }\Omega,&&
\\\label{poro-stat-el2}
&\pl_c^{}\varphi(e,c)+z\phi\ni\mu=\text{\,some constant}&&\text{on }\Omega
,&&
\\\label{poro-stat-el4}
&{\rm div}\big(\eps\nabla\phi)+z\Cdot c=z_{_{\rm DOP}}+{\rm div}(z_{_{\rm DOP}}u)
&&\text{on }\R^d,
\\[-.2em]\label{poro-stat-el5}
&\!\int_\Omega c\,\d x=C_{\rm total}\,\text{=\,a given constant $\in\R^N$.}
\hspace*{-3em}&&&&
\end{align}\end{subequations}
The right-hand side $z_{_{\rm DOP}}\nabla\phi$ in \eq{poro-stat-el1} is the 
Lorenz force acting on a charged elastic solid in the electrostatic 
field. In \eq{poro-stat-el4}, $\eps=\eps(x)>0$ denotes the permittivity
and the equation \eq{poro-stat-el4} itself is the rest of the full Maxwell
system if all evolution and magnetic effects are neglected. Note that 
\eq{poro-stat-el4} is to be solved on the whole universe with the
natural ``boundary'' condition $\phi(\infty)=0$, assuming
naturally that $z$, $c$, and $z_{_{\rm DOP}}$ are extended on 
$\R^d{\setminus}\Omega$ by zero. Actually, the physical units 
are fixed for notational simplicity in such a way that the Faraday constant 
(which should multiply the charges in \eq{poro-stat-el4}) equals 1.

It is interesting that the underlying potential is not convex and, 
instead of a minimizer as in \eq{swelling-small-static-funct-u-c},
we are now to seek a more general critical point, namely a saddle point 
solving the variational problem:
\begin{align}\label{poro-small-static-el}
\hspace*{-4em}\left.\begin{array}{ll}
\text{Min/max}\!\!\!\!&\displaystyle{(u,c,\phi)\mapsto
\int_\Omega\Big(\varphi\big(e(u),c\big)+
(z\Cdot c{-}z_{_{\rm DOP}})\phi+z_{_{\rm DOP}}\nabla\phi\Cdot u}
\\[-.2em]
&\hspace{8.8em}\displaystyle{
-f\Cdot u\Big)\,\d x-\int_{\R^d}\frac\eps2|\nabla\phi|^2\,\d x-\int_{\GNeu}\!\!g\Cdot u\,\d S,
}\!\!
\\\text{subject to}\!\!\!\!&\displaystyle{\int_\Omega\!c\,\d x=C^{}_{\rm total}},
\ u\In H_\Dir^1(\Omega;\R^d),\ c\In L^1(\Omega;\R^N),
\ \phi\In H^1(\R^d).\!\!\!\!
\end{array}\right\}\!\!
\end{align}
This convex/concave structure is sometimes referred under the name of
electrostatic Lagrangian \cite[Sect.3.2]{ProWet09PEMF}
and is consistent with a convex structure of the internal energy, 
cf.\ the argumentation in Remark~\ref{rem-elimination} or 
\eq{poro-therm-engr-}--\eq{poro-therm-engr} below.

\begin{proposition}\label{prop-static-el}
Let 
$\varphi:\R_{\rm sym}^{d\times d}\times\R^N\to\R\cup\{\infty\}$ 
and $\GDir$, $f$, and $g$ be qualified
as in Proposition~\ref{prop-static},
$z_{_{\rm DOP}}\In L^\infty(\Omega)$, $z\In L^\infty(\Omega;\R^N)$, and
$\eps\in L^\infty(\R^d)$ have a positive infimum.
Moreover, let again 
$\alpha=0$
and 
$C^{}_{\rm total}\In\R^N$ be given so that
$\varphi(0,C^{}_{\rm total}/{\rm meas}_d(\Omega))<\infty$.
Then
the boundary-value problem \eq{poro-stat-el} with 
boundary conditions \eq{poro-steady-BC-u}
possesses 
at least one variational solution 
$(u,c)\in H_\Dir^1(\Omega;\R^d)\times L^{1+\epsilon}(\Omega;\R^N)$ whose 
$\R^N$-valued electrochemical potential 
is constant over $\Omega$ in the sense that 
there exists a constant $\bar\mu\In\pl_c^{}\varphi(e(u),c)+z\phi$ 
on $\Omega$ in the variational sense like in Proposition~\ref{prop-static}.
\end{proposition}

\BP
Existence of a saddle point in this problem is to be seen by the classical 
(Banach-space-valued extension \cite{Fan52FPMT} of)
von~Neumann theorem \cite{Neum28TG} (see also \cite[Ch.49]{Zeid86NFAA})
as well as that it yields
some solution to the system \eq{poro-stat-el} with $\mu$ being the 
(vector-valued) Lagrange multiplier to the constraint in 
\eq{poro-small-static-el}. The assumption 
$\varphi(0,C^{}_{\rm total}/{\rm meas}_d(\Omega))<\infty$
again guarantees the feasibility of \eq{poro-small-static-el}.
\EP

\section{General steady-state problems}\label{sec-steady}

The peculiarity behind the general non-static steady-state problem \eq{EL-poro} 
is that it mixes stored energy and the dissipation energy, cf.\ 
\eq{poro-therm-engr-} below. Thus one should not expect a simple 
variational structure which is usual in problems governed merely 
by stored energy. To illustrate this peculiarity more, let us consider 
$\varphi(e,c)=\frac12\bbC(e{-}Ec)\Colon(e{-}Ec)$ and write formally
the underlying operator when ignoring the boundary condition, i.e. 
\begin{align}
\left(\begin{array}{c}\!u\!\\\!c\!\end{array}\right)
\mapsto\left(\!\begin{array}{c}-{\rm div}(\bbC(e(u){-}Ec))
\\
-\Delta(\bbC E\Colon(Ec{-}e(u)))\end{array}\!\right).
\end{align}
This linear operator is obviously nonsymmetric (thus does not have
any potential) and nonmonotone (and even not pseudomonotone) due to the 
3rd-order term $-\Delta\bbC E\Colon e(u)$. 
Therefore, standard methods does not seem to be applicable. Yet, 
advantageously, the variational structure of the static problems in 
Section~\ref{sec-static} can be combined with a carefully constructed fixed 
point. 

\begin{proposition}[Existence of steady states]\label{prop-steady}
Let $\varphi:\R_{\rm sym}^{d\times d}\times\R\to\R\cup\{\infty\}$ be lower 
semicontinuous, strictly convex and coercive in the sense that 
$\varphi(e,c)\ge\epsilon|e|^2+\epsilon|c|^q$ for some $q>{2^*}'$ and 
$\epsilon>0$, 
$\bbM:\R\to\R^{d\times d}$ is continuous, bounded, and uniformly positive definite,
$\eps\in L^\infty(\R^d)$ have a positive infimum,
$\alpha\ge0$ with $\alpha>0$ on a positive-measure part of $\Gamma$, 
${\rm meas}_{d-1}(\GDir)>0$,
$f\In L^{{2^*}'}\!(\Omega;\R^d)$, $g\In L^{{2^\sharp}'}\!\!(\GNeu;\R^d)$, and 
$\mu_{\rm ext}\In L^{{2^\sharp}'}\!\!(\Gamma)$. 
Then
the boundary-value problem \eq{EL-poro} with 
boundary conditions \eq{poro-steady-BC}
possesses at least one variational solution 
$(u,c)\in H_\Dir^1(\Omega;\R^d)\times L^q(\Omega)$ with the corresponding 
chemical potential $\mu\In H^1(\Omega)\cap L^\infty(\Omega)$.
\end{proposition}

\BP
We construct the single-valued mapping $\tilde\mu\mapsto(u,c)\mapsto\mu$ 
for which the Schauder fixed-point theorem will be used. First, fixing 
$\tilde\mu\In H^1(\Omega)$, we solve 
\begin{align}\label{min-u-c}
\hspace*{-1em}\left.\begin{array}{ll}
\text{Minimize }&\displaystyle{(u,c)\mapsto
\int_\Omega\varphi(e(u),c)-\tilde\mu c-f\Cdot u\,\d x
-\int_{\GNeu}\!\!g\Cdot u\,\d S} 
\\[.5em]\text{subject to }&
u\In H_\Dir^1(\Omega;\R^d)\ \ \text{ and }\ \ c\In L^q(\Omega).\ \ 
\end{array}\right\}
\end{align}
Note that the term $\tilde\mu c$ is integrable and 
$(\tilde\mu,c)\mapsto\tilde\mu c: H^1(\Omega)\times L^q(\Omega)\to L^1(\Omega)$ 
is (weak,weak)-continuous due to the condition $q>{2^*}'$ and the Rellich 
theorem. Due to the assumed strict convexity of $\varphi$, this problem has a 
unique
solution $(u,c)$. It is also important that this solution depends continuously 
on $\tilde\mu$ in the sense that $\tilde\mu\mapsto(u,c):H^1(\Omega)
\to H_\Dir^1(\Omega;\R^d)\times L^q(\Omega)$
is (weak,strong)-continuous, which can be seen when exploiting the assumed 
strict convexity of $\varphi$, cf.\ \cite{Visi84SCRR}.

Having $c\In L^1(\Omega)$, we then
\begin{align}\label{min-mu}
\hspace*{-0em}\left.\begin{array}{ll}
\text{Minimize }&\displaystyle{\mu\mapsto
\int_\Omega\frac12\bbM(c)\nabla\mu\Cdot\nabla\mu\,\d x
+\int_{\GNeu}\frac\alpha2\mu^2-\alpha\mu_{\rm ext}^{}\mu\,\d S} 
\\[.4em]\text{subject to }&
\mu\In H^1(\Omega).\ \ 
\end{array}\right\}
\end{align}
Due to the assumed positive definiteness of $\bbM$ and the (partial) 
positivity of $\alpha$, the problem \eq{min-mu} has a unique solution $\mu$. It
is important that the mapping $c\mapsto\mu:L^1(\Omega)\to H^1(\Omega)$
is (strong,weak)-continuous. Actually, even (strong,strong)-continuity can 
easily be proved but it is not needed for our fixed-point argument. 

It should be emphasized that $\mu$ from \eq{min-mu} does not need to be
a chemical potential corresponding to $(u,c)$. Yet, we will show
that it is if $\mu=\tilde\mu$. Such pair $(\mu,\tilde\mu)$ does 
exists due to the Schauder fixed point theorem. Here we also used that 
the solution $\mu$ ranges an a-priori bounded set in $H^1(\Omega)$ because
$\bbM(\cdot)$ is assumed uniformly positive definite.

The 1st-order optimality conditions for \eq{min-u-c} compose from the
partial G\^ateaux derivatives with respect to $u$ and to $c$ to vanish.
The former condition means the Euler-Lagrange 
equation representing the weak formulation of the boundary-value problem: 
\begin{subequations}\label{u-BVP}
\begin{align}
&{\rm div}\,\pl_e\varphi(e(u),c)+f=0&&\text{on }\ \Omega,&&&&&&
\\&u=0\ \text{ on }\GDir\ \ \text{ and }\ \ \sigma\vec{n}=g&&\text{on }\GNeu,
\end{align}
\end{subequations}
while the latter conditions written for $\tilde\mu=\mu$ yields 
\begin{align}
\pl_c\varphi(e(u),c)-\mu\ni0\qquad\text{on }\ \Omega.
\end{align}
The 1st-order optimality conditions for \eq{min-mu} means the Euler-Lagrange 
equation representing the weak formulation of the boundary-value problem: 
\begin{subequations}\label{mu-BVP}\begin{align}
&{\rm div}(\bbM(c)\nabla\mu)=0&&\text{on }\ \Omega,&&&&&&
\\&
\bbM(c)\nabla\mu\cdot\vec{n}+\alpha\mu=\alpha\mu_{\rm ext}&&\text{on }\ \Gamma.
\end{align}
\end{subequations}
Altogether, \eq{u-BVP}--\eq{mu-BVP} reveal that $(u,c)$ solves 
the 
problem 
\eq{EL-poro}--\eq{poro-steady-BC}.
\EP

Let us now investigate the steady-state variant of the static 
electrically-charged multi-component problem \eq{poro-stat-el}. This results to
\begin{subequations}\label{poro-steady-el}
\begin{align}\label{poro-steady-el1}
&{\rm div}\,\pl_e\varphi(e(u),c)+f=z_{_{\rm DOP}}\nabla\phi
&&\text{on }\Omega,&&
\\\label{poro-steady-el2}
&{\rm div}(\bbM(c)\nabla\mu)+r(c)=0\ \ \ \text{with }\ 
\mu\in\pl_c\varphi(e(u),c)+z\phi&&\text{on }\Omega
,&&
\\\label{poro-steady-el4}
&{\rm div}\big(\eps\nabla\phi)+z\Cdot c=z_{_{\rm DOP}}+{\rm div}(z_{_{\rm DOP}}u)
&&\text{on }\R^d,
\hspace*{-3em}&&&&
\end{align}\end{subequations}
with $r=r(c)$ the rate of chemical reactions, to be completed by the boundary 
conditions \eq{poro-steady-BC} and $\phi(\infty)=0$.

\begin{proposition}\label{prop-steady-el}
Let again $\varphi:\R_{\rm sym}^{d\times d}\times\R^N\to\R\cup\{\infty\}$ be lower
semicontinuous, strictly convex, and now even uniformly convex in $c$ and 
coercive in the sense that, for some 
$\epsilon\!>\!0$, 
\begin{subequations}\begin{align}\nonumber
&\hspace*{-4em}\forall\,(e_1,c_1),(e_2,c_2)\In\Rsym\times\R^N\ \ \forall\,
m_1\In\pl_c^{}\varphi(e_1,c_1),\ m_2\In\pl_c^{}\varphi(e_2,c_2):
\\&\hspace*{-4em}\nonumber
\ \ \ \epsilon
\big(|c_1|^{q-2}c_1-|c_2|^{q-2}c_2\big)\big(c_1{-}c_2\big)\le
\big(\pl_e^{}\varphi(e_1,c_1)-\pl_e^{}\varphi(e_2,c_2)\big)
\Colon (e_1{-}e_2)
\\[-.1em]&\hspace*{-4em}\hspace*{21.7em}
+
(m_1{-}m_2)\Cdot(c_1{-}c_2),
\label{poro-unif-convex-of-phi}
\\&\hspace*{-4em}\exists\,q\!>\!{2^*}'\ \:
\forall\,(e,c)\In\Rsym\times\R^N:\qquad
\varphi(e,c)\ge\epsilon|e|^2+\epsilon|c|^q.
\end{align}\end{subequations}
Let further $\bbM:\R^N\to\R^{d\times d\times N}$ be continuous, bounded, and 
uniformly positive definite,
$r:\R^N\to\R^N$ continuous and bounded,
$\eps\in L^\infty(\R^d)$ have a positive infimum,
$\alpha\ge0$ with $\alpha>0$ on a positive-measure part of $\Gamma$, 
$z_{_{\rm DOP}}\In L^\infty(\Omega)$, $z\In L^\infty(\Omega;\R^N)$, 
${\rm meas}_{d-1}(\GDir)>0$,
$f\In L^{{2^*}'}\!(\Omega;\R^d)$, $g\In L^{{2^\sharp}'}\!\!(\GNeu;\R^d)$, and 
$\mu_{\rm ext}\In L^{{2^\sharp}'}\!\!(\Omega;\R^N)$. Then
the boundary-value problem \eq{poro-steady-el} with 
boundary conditions \eq{poro-steady-BC} and $\phi(\infty)=0$ possesses 
at least one variational solution 
$(u,c)\in H_\Dir^1(\Omega;\R^d)\times L^q(\Omega;\R^N)$
with the corresponding electrochemical potential $\mu\In H^1(\Omega;\R^N)$.
\end{proposition}

\BP
We construct the single-valued mapping $\tilde\mu\mapsto(u,c,\phi)\mapsto\mu$ 
for which the Schauder fixed-point theorem will be used. First, we fix 
$\tilde\mu\In H^1(\Omega;\R^N)$ and, being motivated by 
\eq{poro-small-static-el}, we modify \eq{min-u-c} as
\begin{align}\label{min/max-u-c-phi}
\hspace*{-4em}\left.\begin{array}{ll}
\text{Min/max}&\displaystyle{(u,c,\phi)\mapsto
\int_\Omega\Big(\varphi(e(u),c)+
(z\Cdot c{-}z_{_{\rm DOP}})\phi+z_{_{\rm DOP}}\nabla\phi\Cdot u}
\\[-.4em]
&\hspace{5.em}\displaystyle{-\tilde\mu\Cdot c-f\Cdot u\Big)\,\d x
-\int_{\R^d}\frac\eps2|\nabla\phi|^2\,\d x-\int_{\GNeu}\!\!g\Cdot u\,\d S,}
\\\text{subject to}&
u\In H_\Dir^1(\Omega;\R^d),\ \ c\In L^q(\Omega;\R^N),
\ \ \phi\In H^1(\R^d).\ \ 
\end{array}\right\}\!\!
\end{align}
Due to the assumed strict convexity of $\varphi$, this problem has a unique
solution $(u,c,\phi)$ which depends continuously on $\tilde\mu$ in the sense 
that 
$\tilde\mu\mapsto(u,c,\phi):H^1(\Omega;\R^N)\to H_\Dir^1(\Omega;\R^d)\times 
L^q(\Omega;\R^N)\times H^1(\R^d)$ is (weak,strong)-continuous. More in 
detail, the existence of this saddle point is the classical 
(Banach-space-valued extension \cite{Fan52FPMT} of)
von~Neumann theorem \cite{Neum28TG}.
The uniqueness can be proved by analyzing the optimality conditions
\eq{poro-steady-el1} and \eq{poro-steady-el4} together with 
$\pl_c^{}\varphi(e(u),c)+z\phi\ni\tilde\mu$ written as
\begin{align}\label{chem-pot-tilde}
m+z\phi=\tilde\mu\qquad\text{for some}\ \ m\in\pl_c^{}\varphi(e(u),c),
\end{align}
considered for two solutions $(u_i,c_i,\phi_i,m_i)$, $i=1,2$, and subtracted.
Using the abbreviation $u_{12}=u_1{-}u_2$, $c_{12}=c_1{-}c_2$, 
$\phi_{12}=\phi_1{-}\phi_2$,  and $m_{12}=m_1{-}m_2$, this results to the system
\begin{subequations}\label{BVP-unique}\begin{align}\label{BVP-unique-u}
&{\rm div}\big(\pl_e\varphi(e(u_1),c_1)-\pl_c^{}\varphi(e(u_2),c_2)\big)
=z_{_{\rm DOP}}\nabla\phi_{12}
&&\text{on }\Omega,&&
\\\label{BVP-unique-phi}
&{\rm div}\big(\eps\nabla\phi_{12})+z\Cdot c_{12}={\rm div}(z_{_{\rm DOP}}u_{12})
&&\text{on }\R^d,
\\\label{BVP-unique-mu}
&
m_{12}+z\phi_{12}=0&&\text{on }\Omega,&&
\end{align}\end{subequations}
with the homogeneous boundary conditions for \eq{BVP-unique-u},
i.e.\ $u_{12}=0$ on $\GDir$ and 
$(\pl_e\varphi(e(u_1),c_1)-\pl_e^{}\varphi(e(u_2),c_2))\vec{n}=0$ on $\GNeu$.
Testing \eq{BVP-unique-u} by $u_{12}$ and also \eq{BVP-unique-mu} by $c_{12}$ 
(integrated it over $\Omega$)
and \eq{BVP-unique-phi} by $\phi_{12}$  (integrated it over $\R^d$) gives 
\begin{align}\nonumber
&\hspace*{-4em}\int_\Omega\!\big(\pl_e^{}\varphi(e(u_1),c_1)-\pl_e^{}\varphi(e(u_2),c_2)\big)
\Colon e(u_{12})+m_{12}\Cdot c_{12}\,\d x
+\int_{\R^d}\!\!\eps|\nabla\phi_{12}|^2\,\d x
\\[-.7em]\nonumber&\hspace*{-4em}\quad=\int_\Omega\! 
-z_{_{\rm DOP}}\nabla\phi_{12}\Cdot u_{12}-z\Cdot c_{12}\,\phi_{12}\,\d x
+\int_{\R^d}\!\!z\Cdot c_{12}
\phi_{12}-{\rm div}(z_{_{\rm DOP}}u_{12})\phi_{12}\,\d x
\\[-.7em]\nonumber&\hspace*{-4em}\quad\qquad=\int_{\R^d}\!\!
-z_{_{\rm DOP}}\nabla\phi_{12}\Cdot u_{12}-z\Cdot c_{12}\,\phi_{12}
+z\Cdot c_{12}\phi_{12}-{\rm div}(z_{_{\rm DOP}}u_{12})\phi_{12}\,\d x
\\[-.7em]&\hspace*{-4em}\qquad\qquad=\int_{\R^d}\!\!-{\rm div}(z_{_{\rm DOP}}u_{12}\phi_{12})\,\d x
=-\big[z_{_{\rm DOP}}u_{12}\phi_{12}\big](\infty)=0.
\label{poro-est-for-uniqueness}\end{align}
We used that $z=0$ on $\R^d{\setminus}\Omega$ and then cancellation of the 
terms $\pm z\Cdot c_{12}\phi_{12}$ as well as that 
$\phi_{12}(\infty)=0$ and $z_{_{\rm DOP}}(\infty)=0$. 

From the strict monotonicity of 
$\pl\varphi$, we can easily see uniqueness of 
the saddle point of \eq{min/max-u-c-phi},
needed for the Schauder fixed point. Moreover, the mentioned continuity 
is to be proved by taking two right-hand sides $\tilde\mu_i$, $i=1,2$, in 
\eq{chem-pot-tilde}. This gives rise the additional term 
$\int_\Omega\tilde\mu_{12}\Cdot c_{12}\d x$ on the right-hand side of 
\eq{poro-est-for-uniqueness}, which can be estimated by using the
H\"older inequality and then, by the assumption 
\eq{poro-unif-convex-of-phi} and the uniform convexity of the 
$L^q(\Omega;\R^N)$-space, again obtain the desired (weak,strong)-continuity
of $\tilde\mu\mapsto c:H^1(\Omega;\R^N)\to L^q(\Omega;\R^N)$; cf.\ 
e.g.\ \cite[Chap.\,2]{Roub13NPDE}. 

Having $c\In L^1(\Omega;\R^N)$, we again solve \eq{min-mu} now 
in addition with a term $-r(c)\Cdot\mu$ and, as in 
the proof of Proposition~\ref{prop-steady}, we prove that 
the mapping $c\mapsto\mu:L^1(\Omega;\R^N)\to H^1(\Omega;\R^N)$
is (strong,weak)-continuous. 

It should be emphasized that $\mu$ from \eq{min-mu} does not need to be
a chemical potential corresponding to $(u,c)$. Yet, it is if $\mu=\tilde\mu$. 
Such pair $(\mu,\tilde\mu)$ does 
exist due to the Schauder fixed point theorem. Here we also used that 
the solution $\mu$ ranges in an a-priori bounded subset of $H^1(\Omega)$ because
of the assumed uniform positive definiteness of $\bbM(\cdot)$ and the 
boundedness of the reaction rates $r(\cdot)$. 
\EP

\section{Anisothermal problems}\label{sec-temp}

We already mentioned that the diffusion equation \eq{EL-poro-mu} or
\eq{poro-steady-el2} is 
related to the dissipation rather than the stored energy. Thermodynamically,
the dissipation rate (i.e.\ here $\bbM\nabla\mu\Cdot\nabla\mu$)
leads to the heat production which might substantially influence 
temperature if the specimen is large 
or/and the produced heat cannot be transferred away sufficiently fast. 
In turn, variation of temperature may influence the dissipation mechanism and 
the stored energy too, and thus gives rise to a 
{\it thermomechanically coupled system}.

The free energy $\varphi$ as well as the mobility tensor $\bbM$ 
now may depend on temperature, let us denote it by $\theta$. The original 
system \eq{EL-poro} then augments as
\begin{subequations}\label{EL-poro-thermo}\begin{align}
&{\rm div}\,\sigma+f=0&&\text{with }\ \ \ \sigma=\pl_e\varphi(e(u),c,\theta),&&
\\\label{EL-poro-thermo2}
&{\rm div}(\bbM(c,\theta)\nabla\mu)=0&&\text{with }\ \ \ 
\mu\in\pl_c\varphi(e(u),c,\theta),&&
\\\label{EL-poro-thermo3}
&{\rm div}(\bbK(c,\theta)\nabla\theta)+\bbM(c,\theta)\nabla\mu\Cdot\nabla\mu
=0\hspace*{-5em}
\end{align}
\end{subequations}
to be solved on a bounded Lipschitz domain $\Omega\subset\R^d$.
Note that \eq{EL-poro-thermo3} involves the {\it Fourier law}, saying
that the heat flux equals $-\bbK(c,\theta)\nabla\theta$.
This system should be completed by suitable boundary conditions, e.g.\ 
\begin{subequations}\label{poro-steady-BC-thermo}
\begin{align}\label{poro-therm-BC-u}
&u=0\ \text{ on }\GDir,\qquad\ \ \sigma\vec{n}=g\ \ \text{ on }\GNeu,
\quad\ 
\\&\label{poro-therm-BC-c}
\bbM(c,\theta)\nabla\mu\cdot\vec{n}+\alpha\mu=\alpha\mu_{\rm ext}\ \text{ on }\Gamma,
\\&\label{poro-therm-BC-theta}
\bbK(c,\theta)\nabla\theta\cdot\vec{n}
+\gamma\theta=\gamma\theta_{\rm ext}\ \ \ \text{ on }\Gamma.
\end{align}
\end{subequations}
In \eq{EL-poro-thermo3} and \eq{poro-therm-BC-theta},
$\bbK=\bbK(c,\theta)$ denotes a heat-conductivity tensor.

In this scalar case, an interesting transformation (used also
in a steady-state thermistor problem \cite[Sect.\,6.4]{Roub13NPDE})
is based on the formula 
${\rm div}(av)=a\,{\rm div}\,v+\nabla a\cdot v$. One can indeed rely on 
\begin{align}
\hspace*{-4em}\!\!{\rm div}\big(\mu\bbM(c,\theta)\nabla\mu\big)=
\mu\!\!\!\!\!\lineunder{{\rm div}\big(\bbM(c,\theta)\nabla\mu\big)}{$=0$ by 
\eq{EL-poro-thermo2}}\!\!\!\!\!\!
+\bbM(c,\theta)\nabla\mu\Cdot\nabla\mu=\bbM(c,\theta)\nabla\mu\Cdot\nabla\mu.\
\label{calculus}\end{align}
If $\mu\in L^\infty(\Omega)$, it indeed can be tested by functions from 
$H^1(\Omega)$ and thus lives in $H^1(\Omega)^*$. In the scalar case, one
has the information $\mu\in L^\infty(\Omega)$ at disposal due to the 
maximum principle if the external chemical potential $\mu_{\rm ext}$ is
in  $L^\infty(\Gamma)$. Thus, instead of \eq{EL-poro-thermo3},
one can equivalently consider 
\begin{align}\label{transformed-heat-eq}
{\rm div}\big(\bbK(c,\theta)\nabla\theta+\mu\bbM(c,\theta)\nabla\mu\big)
=0.
\end{align}

\begin{proposition}[Existence of steady states]\label{prop-steady-therm}
Let $\varphi:\R_{\rm sym}^{d\times d}\times\R\times\R\to\R\cup\{\infty\}$ be 
lower semicontinuous and coercive in the sense that 
$\varphi(e,c,\theta)\ge\epsilon|e|^2+\epsilon|c|^q$ for some 
$q>{2^*}'$ and $\epsilon>0$ and $\varphi(\cdot,\cdot,\theta)$ be strictly 
convex for any $\theta\in\R$ and 
$\theta\mapsto\varphi(\cdot,\cdot,\theta)$ be continuous in the sense
of $\varGamma$-convergence as stated in \eqref{Gamma-convergence} below,
and let $\bbM,\bbK:\R^2\to\R^{d\times d}$ be continuous, bounded, and uniformly 
positive definite, $\gamma\ge0$, $\mu_{\rm ext}\In L^\infty(\Gamma)$,
$\theta_{\rm ext}\In L^{{2^\sharp}'}\!\!(\Gamma)$,
and $\alpha$, $f$, and $g$ be as in Proposition\,\ref{prop-steady}.
Then
the boundary-value problem \eq{EL-poro-thermo}--\eq{poro-steady-BC-thermo}
possesses 
at least one variational solution 
$(u,c,\theta)\in H_\Dir^1(\Omega;\R^d)\times L^q(\Omega)\times H^1(\Omega)$ with 
the corresponding chemical potential $\mu\In H^1(\Omega)\cap L^\infty(\Omega)$. 
\end{proposition}

\BP
We construct the single-valued mapping 
$(\tilde\mu,\tilde\theta)\mapsto(u,c,\theta)\mapsto\mu$ for which the Schauder 
fixed-point theorem will be used. First, fixing 
$\tilde\mu\In H^1(\Omega)\cap L^\infty(\Omega)$ and 
$\tilde\theta\In H^1(\Omega)$, we solve 
\begin{align}\label{min-u-c-therm}
\hspace*{-4em}\left.\begin{array}{ll}
\text{Minimize }&\displaystyle{(u,c)\mapsto
\int_\Omega\varphi(e(u),c,\tilde\theta)
-\tilde\mu c-f\Cdot u\,\d x
-\int_{\GNeu}\!\!g\Cdot u\,\d S} 
\\[.9em]\text{subject to }&
u\In H_\Dir^1(\Omega;\R^d)\ \ \text{ and }\ \ c\In L^q(\Omega).\ \ 
\end{array}\right\}
\end{align}
Due to the assumed strict convexity of $\varphi$, this problem has a unique
solution $(u,c)$. It is also important that this solution depends continuously 
on $(\tilde\mu,\tilde\theta)$ in the sense that 
$(\tilde\mu,\tilde\theta)\mapsto(u,c):
H^1(\Omega)^2\to H_\Dir^1(\Omega;\R^d)\times L^q(\Omega)$
is (weak,strong)-continuous. In particular, we use 
the assumed $\Gamma$-convergence, meaning that the set-valued mapping
\begin{align}\label{Gamma-convergence}
\hspace*{-2em}\theta\mapsto{\rm epi}\,\varphi(\cdot,\cdot,\theta)
:=\big\{(e,c,a)\In \R_{\rm sym}^{d\times d}{\times}\R^N{\times}(\R\cup\{\infty\});\ 
\varphi(e,c,\theta)\le a\big\}
\end{align}
is continuous in the Hausdorff sense, and also the Rellich compactness
theorem so that the functional 
$$
H_\Dir^1(\Omega;\R^d)\times L^q(\Omega)\to
\R\cup\{\infty\}\ :\ (u,c)\mapsto\int_\Omega\!\varphi(e(u),c,\tilde\theta)
-\tilde\mu c\,\d x
$$ 
$\Gamma$-converges if $(\tilde\mu,\tilde\theta)$ converges weakly in 
$H^1(\Omega)^2$. From this and the strict convexity, the desired 
continuity of  $(\tilde\mu,\tilde\theta)\mapsto(u,c)$ is seen.

Further, we solve the boundary-value problem \eq{mu-BVP} now with 
$\bbM(c,\tilde\theta)$ instead of $\bbM(c)$. In addition, assuming 
$\mu_{\rm ext}\In L^\infty(\Gamma)$, we can use the maximum principle yielding 
the estimate ${\rm ess\,inf}\,\mu_{\rm ext}(\Gamma)\le\mu(x)
\le{\rm ess\,sup}\,\mu_{\rm ext}(\Gamma)$ for a.a.\ $x\In\Omega$.
This estimate is independent of $(\tilde\mu,\tilde\theta)$,
as well as the estimates $\|\mu\|_{H^1(\Omega)}\le C$ and 
$\|\theta\|_{H^1(\Omega)}\le C$ provided $C$ is large enough. 

Eventually, having  $c$, $\mu$, and $\tilde\theta$ at disposal, we solve
\begin{align}\label{min-theta}
\left.\begin{array}{ll}
\text{Minimize }&\displaystyle{\theta\mapsto
\int_\Omega\!\Big(\frac12\bbK(c,\tilde\theta)\nabla\theta
+\mu\bbM(c,\tilde\theta)\nabla\mu\Big)\Cdot\nabla\theta\,\d x}
\\[-.5em]&\displaystyle{\hspace*{9.5em}
+\int_\Gamma\frac\gamma2\theta^2-\gamma\theta_{\rm ext}\theta\,\d S}
\\\text{subject to }&
\theta\In H^1(\Omega).
\end{array}\right\}
\end{align}
For the (strong$\times$strong$\times$weak,weak)-continuity of the mapping 
$(c,\tilde\theta,\tilde\mu)\mapsto\theta$, it 
is important that the weak convergence of $\mu$ in $H^1(\Omega)$ 
implies that ${\rm div}(\mu\bbM(c,\tilde\theta)\nabla\mu)$
converges weakly in $H^1(\Omega)^*$ so that the weak convergence 
temperatures in $H^1(\Omega)$ easily follows. 

This allows to execute the Schauder fixed point for the mapping 
$(\tilde\mu,\tilde\theta)\mapsto(\mu,\theta):H^1(\Omega)^2\to H^1(\Omega)^2$
in the weak topology. (In fact, even strong convergence of $\mu$'s solving 
\eq{mu-BVP} now with $\bbM=\bbM(c,\tilde\theta)$ and thus also of $\theta$'s
solving \eq{transformed-heat-eq} with $\bbK=\bbK(c,\tilde\theta)$
and $\bbM=\bbM(c,\tilde\theta)$ can be proved but, in contrast with the
proof of Proposition~\ref{prop-steady-el-therm}, we will not need it here.)
\EP

The thermodynamical completion of the electrically-charged multicomponent 
system  combines \eq{poro-steady-el} with 
\eq{EL-poro-thermo}, resulting to 
\begin{subequations}\label{poro-steady-el-therm}
\begin{align}\label{poro-steady-el-therm1}
&\hspace*{-2em}{\rm div}\,\pl_e\varphi(e(u),c,\theta)+f=z_{_{\rm DOP}}\nabla\phi
&&\text{on }\Omega,&&
\\\label{poro-steady-el-therm2}
&\hspace*{-2em}{\rm div}(\bbM(c,\theta)\nabla\mu)+r(c,\theta)=0\ \ \ \text{with }\ 
\mu\in\pl_c\varphi(e(u),c,\theta)+z\phi&&\text{on }\Omega
,&&
\\\label{poro-steady-el-therm3}
&\hspace*{-2em}{\rm div}(\bbK(c,\theta)\nabla\theta)+\bbM(c,\theta)\nabla\mu\Colon\nabla\mu
+h(c,\theta)
=\mu\Cdot r(c,\theta)&&\text{on }\Omega,\hspace*{-5em}
\\\label{poro-steady-el-therm4}
&\hspace*{-2em}{\rm div}\big(\eps\nabla\phi)+z\Cdot c=z_{_{\rm DOP}}+{\rm div}(z_{_{\rm DOP}}u)
&&\text{on }\R^d.
\hspace*{-3em}&&&&
\end{align}\end{subequations}
The right-hand side $\mu\Cdot r(c,\theta)$ of 
\eq{poro-steady-el-therm3} represents the (negative) heat production
where $h=h(c,\theta)$ denotes the heat-production rate 
due to chemical reactions. 
Exploiting \eq{poro-steady-el-therm2}, the calculus \eq{calculus} modifies to  
${\rm div}\big(\mu\bbM(c,\theta)\nabla\mu\big)=
\bbM(c,\theta)\nabla\mu\Colon\nabla\mu-\mu\Cdot r(c,\theta)$
so that the heat-transfer problem \eq{poro-steady-el-therm3} turns again to 
\eq{transformed-heat-eq}. The multi-component system \eq{poro-steady-el-therm} 
is more complicated than \eq{EL-poro-thermo} because 
the maximum principle for $\mu$ and the variational structure for the
heat-transfer equation (even if transformed into \eq{transformed-heat-eq})
is not at disposal, however. 

\begin{proposition}\label{prop-steady-el-therm}
Let $\varphi:\R_{\rm sym}^{d\times d}\times\R^N\times\R\to\R\cup\{\infty\}$ be 
lower semicontinuous and coercive in the sense that 
$\varphi(e,c,\theta)\ge\epsilon|e|^2+\epsilon|c|^q$ for some 
$q>{2^*}'$ and $\epsilon>0$ and $\varphi(\cdot,\cdot,\theta)$ be strictly 
convex satisfying \eq{poro-unif-convex-of-phi} uniformly for 
any $\theta\in\R$ and $\theta\mapsto\varphi(\cdot,\cdot,\theta)$
be continuous in the sense of $\varGamma$-convergence as stated previously in
\eqref{Gamma-convergence},
and let $\bbM:\R^N\times\R\to\R^{d\times d\times N}$
and $\bbK:\R^N\times\R\to\R^{d\times d}$ be continuous, bounded, and uniformly 
positive definite, $h:\R^N\times\R\to\R$ be continuous and bounded,
$\gamma\ge0$, and
$\theta_{\rm ext}\In L^1
(\Gamma)$. Moreover, let $\mu_{\rm ext}$, $\eps$, $\alpha$, $z_{_{\rm DOP}}$, 
$z$, $f$, and $g$ be as in Proposition~\ref{prop-steady-el}.
Then
the boundary-value problem \eq{poro-steady-el-therm} with 
boundary conditions \eq{poro-steady-BC} and $\phi(\infty)=0$
possesses 
at least one variational solution 
$(u,c,\theta,\phi)\In H_\Dir^1(\Omega;\R^d)\times L^q(\Omega;\R^N)\times 
W^{1,p}(\Omega)\times H^1(\R^d)$ with any $1\le p<
d'$ with 
the corresponding electro\-chemical potential $\mu\In H^1(\Omega;\R^N)$. 
\end{proposition}

{\noindent{\it Sketch of the proof.\ }}
We organize the Schauder fixed point for a composed single-valued mapping
\begin{subequations}\label{poro-steady-el-therm-FP}
\begin{align}\label{poro-steady-el-therm-FP1}
&\hspace*{-4em}\!\!\!\!(\tilde\mu,\tilde\theta){\mapsto}(u,c,\phi){:}
H^1(\Omega;\R^N){\times}W^{1,p}(\Omega){\to}
H_\Dir^1(\Omega;\R^d){\times}L^q(\Omega;\R^N){\times}H^1(\R^d),\!\! 
\\\label{poro-steady-el-therm-FP2}
&\hspace*{-4em}\!\!\!\!(c,\tilde\theta)\mapsto\mu:L^q(\Omega;\R^N)\times W^{1,p}(\Omega)\to H^1(\Omega;\R^N),\text{ and eventually}
\\\label{poro-steady-el-therm-FP3}
&\hspace*{-4em}\!\!\!\!(c,\tilde\theta,\mu)\mapsto\theta:L^q(\Omega;\R^N)\times W^{1,p}(\Omega)\times H^1(\Omega;\R^N)\to W^{1,p}(\Omega).
\end{align}\end{subequations}

For \eq{poro-steady-el-therm-FP1}, the minimization variational problem 
\eq{min-u-c-therm} in the previous proof is to be replaced by the 
saddle-point problem \eq{min/max-u-c-phi} but now with 
$\varphi=\varphi(\cdot,\cdot,\tilde\theta)$. 
The uniqueness of its solution is again due to 
\eq{BVP-unique}--\eq{poro-est-for-uniqueness} but now with 
$\varphi(e(u_i),c_i,\tilde\theta)$ instead of $\varphi(e(u_i),c_i)$,
as well as the (weak,weak$\times$strong$\times$weak)-continuity of 
\eq{poro-steady-el-therm-FP1}. 

Moreover, for \eq{poro-steady-el-therm-FP2}, we solve again the 
minimization problem \eq{min-mu} with $\bbM=\bbM(c,\tilde\mu)$
and additionally with the term $\mu\Cdot r(c,\tilde\theta)$.
In contrast to the proofs of 
Propositions~\ref{prop-steady}--\ref{prop-steady-therm}, we now need the 
(strong$\times$weak,strong)-continuity of \eq{poro-steady-el-therm-FP2},
which follows standardly by the uniform convexity of the functional in
\eq{min-mu}.

Eventually, for \eq{poro-steady-el-therm-FP3},
instead of the minimization problem \eq{min-theta}
whose infimum might be $-\infty$ because now
$\tilde\mu\bbM(c,\tilde\theta)\nabla\tilde\mu\Cdot\nabla\theta\not\in L^1(\Omega)$ in general since the $L^\infty$-estimate on $\tilde\mu$ (and thus on 
$\tilde\mu$ too) is not at disposal, we should solve \eq{poro-steady-el-therm3}
with the boundary condition \eq{poro-therm-BC-theta} with 
$\bbK=\bbK(c,\tilde\theta)$, $\bbM=\bbM(c,\tilde\theta)$, and 
the heat sources $h(c,\tilde\theta)+\mu\Cdot r(c,\tilde\theta)$, by the 
non-variational
method. By the classical Stampacchia \cite{Stam65PDEE} transposition method, 
see also e.g.\ \cite[Section 3.2.5]{Roub13NPDE},
this linear boundary-value problem has a unique variational solution 
$\theta$ which belongs to $W^{1,p}(\Omega)$ with any $1\le p<d'=d/(d{-}1)$.
The (strong$\times$weak$\times$strong,weak)-continuity of 
\eq{poro-steady-el-therm-FP3} is obvious.

Altogether, the Schauder fixed-point relies on the weak continuity of the 
mapping $(\tilde\mu,\tilde\theta)\mapsto(\mu,\theta):
H^1(\Omega;\R^N)\times W^{1,p}(\Omega)\to H^1(\Omega;\R^N)\times W^{1,p}(\Omega)$.
\EP

\section{Concluding remarks}\label{sec-rem}

\begin{remark}[{\sl Uniqueness}]\label{rem-uniqueness}
\upshape
Interestingly, it is not clear whether the solution of the static problem 
in Proposition~\ref{prop-static} is unique, even if $\varphi$ is 
strictly convex so that the solution obtained in the proof of 
Proposition~\ref{prop-static} by solving 
\eq{swelling-small-static-funct-u-c} is unique. 
An analogous comment is relevant to the problem in 
Proposition~\ref{prop-static-el}. As far as all steady-state but non-static
problems, the nonuniqueness is rather to be expected.
\end{remark}

\begin{remark}[{\sl About the fixed-point strategy}]\label{rem-other-FP}
\upshape
Reorganizing the fixed point as $\tilde c\mapsto\mu\mapsto(u,c)$ or 
$\tilde c\mapsto\mu\mapsto(u,c,\phi)$ and then
pursuing a ``construction'' of $c=\tilde c$ would 
bring troubles with compactness needed in the Schauder 
fixed-point theorem because any information about $\nabla c$ is missing.
Note that, from \eq{swelling-small-c=c(q)}, we could think about $\nabla c=
\pl_{\mu e}^2\varphi^*(e(u),\mu)\nabla e(u)+\pl_{\mu\mu}^2\varphi^*(e(u),\mu)\nabla\mu$ but it would need to have $\nabla e(u)$ estimated. Therefore, it does not
seem much freedom to organize the proof of the ``non-static'' 
Propositions~\ref{prop-steady}--\ref{prop-steady-el-therm}.
A certain freedom is in \eq{min-theta} which may involve $\tilde\mu$
instead of $\mu$ because the weak convergence of $\mu$'s is enough
for \eq{min-theta}. This freedom is not applicable to \eq{poro-steady-el-therm3}
used in the proof of Proposition~\ref{prop-steady-el-therm}, however.
\end{remark}

\begin{remark}[{\sl Elimination of the saddle-point structure}]\label{rem-elimination}
\upshape
Making maximization with respect to $\phi$ in \eq{poro-small-static-el}
or in \eq{min/max-u-c-phi} may eliminate the $\phi$ variable, cf.\ also  
e.g.\ \cite[Sect.49.2]{Zeid86NFAA} for a general viewpoint.
Using \eq{poro-stat-el4} together with the calculus 
\begin{align}\nonumber
\hspace*{-3em}\int_\Omega z\phi
\Cdot c-z_{_{\rm DOP}}\phi+z_{_{\rm DOP}}\nabla\phi\Cdot u\,\d x
&=\int_{\R^d}\phi\Big(z\Cdot c-z_{_{\rm DOP}}
-{\rm div}(z_{_{\rm DOP}}u)\Big)\,\d x
\\[-.3em]&\label{calculus+}
\hspace*{-3em}=\ -\!\int_{\R^d}\phi\,{\rm div}\big(\eps\nabla\phi)
=\int_{\R^d}\eps|\nabla\phi|^2\,\d x\,,
\end{align}
the convex/concave problem \eq{min/max-u-c-phi} turns into the
convex constrained problem:
\begin{align}\label{min-u-c-no-phi}
\hspace*{-3em}\left.\begin{array}{ll}
\text{Minimize}&\displaystyle{(u,c,\phi)\mapsto
\int_\Omega\varphi(e(u),c)
-\tilde\mu\Cdot c-f\Cdot u\,\d x} 
\\[-.2em]
&\hspace{10em}\displaystyle{
+\int_{\R^d}\frac\eps2|\nabla\phi|^2\,\d x-\int_{\GNeu}\!\!g\Cdot u\,\d S,}
\\[.9em]\text{subject to}&
{\rm div}\big(\eps\nabla\phi)+z\Cdot c=z_{_{\rm DOP}}+{\rm div}(z_{_{\rm DOP}}u)
\ \ \text{on }\R^d,
\\[.5em]&
u\In H_\Dir^1(\Omega;\R^d),\ \ c\In L^1(\Omega;\R^N),\ \ \phi\In H^1(\R^d)
.\ \ 
\end{array}\right\}
\end{align}
\end{remark}

\begin{remark}[{\sl Towards electroneutrality}]
\upshape
We can further eliminate the electrostatic 
potential by introducing the electrical induction $\vecd=\eps\nabla\phi$
and, defining the Banach space 
\begin{align*}
L_{{\rm rot},\eps}^2(\R^d;\R^d):=\big\{&\vecd\In L^2(\R^d;\R^d);\ \:
\exists\,\phi\In H^1(\R^d):\ 
\\&\vecd=\eps\nabla\phi
\text{ in the sense of distributions}\big\},
\end{align*}
rewrite \eq{min-u-c-no-phi} as
\begin{align}\label{min-u-c-no-phi+}
\hspace*{-3em}\left.\begin{array}{ll}
\text{Minimize }&\displaystyle{(u,c,\phi)\mapsto
\int_\Omega\varphi(e(u),c)
-\tilde\mu\Cdot c-f\Cdot u\,\d x} 
\\[-.3em]
&\hspace{11em}\displaystyle{
+\int_{\R^d}\frac1{2\eps}|\vecd|^2\,\d x-\int_{\GNeu}\!\!g\Cdot u\,\d S,}
\\[.3em]\text{subject to }&
{\rm div}\,\vecd+z\Cdot c=z_{_{\rm DOP}}\!+{\rm div}(z_{_{\rm DOP}}u)
\ \ \text{on }\R^d,
\\[.5em]&
u\In H_\Dir^1(\Omega;\R^d),\ \ \ \ c\In L^1(\Omega;\R^N),\ \ \ \ 
\vecd\In L_{{\rm rot},\eps}^2(\R^d;\R^d).\ \ 
\end{array}\right\}
\end{align}
It reveals an asymptotics for $\eps\to0$: namely, assuming 
$\eps(x)=\epsilon\eps_0(x)$,
the space $L_{{\rm rot},\eps}^2(\R^d;\R^d)$
is independent of $\epsilon$ and then 
$\|\vecd\|_{L^2(\R^d;\R^d)}=\mathscr{O}(\epsilon^{1/2})$ for $\epsilon\to0$ and thus 
$$
\big\|z\Cdot c-z_{_{\rm DOP}}\!-{\rm div}(z_{_{\rm DOP}}u)\big\|_{H^{-1}(\R^d;\R^d)}
=\mathscr{O}(\epsilon^{1/2}).
$$ 
In particular, in the limit, one may expect the {\it electroneutrality}, i.e.\ 
$z\Cdot c-z_{_{\rm DOP}}\!-{\rm div}(z_{_{\rm DOP}}u)=0$. Without any
rigorous justification, this ansatz is indeed 
often used in computational implementation if the specimen size is 
substantially bigger than the so-called Debye length to avoid
spatially extremely stiff problems arising for small permittivities $\eps$,
cf.\ e.g.\ \cite{Fuhr??MNMF}.
\end{remark}

\begin{remark}[{\sl Thermodynamics}]
\upshape
Inspecting evolution variant of the anisothermal electrically-charged
problem is illustrative. For simplicity, let us consider the special free 
energy suppressing thermo-mechanical interactions by assuming 
\begin{align}
\psi(e,c,\theta,\vece)
=\psi_{_{\rm ME}}(e,c)+\psi_{_{\rm TH}}(\theta)-\frac\eps2|\vece|^2 
\end{align}
with $\vece$ standing for the intensity of the electrostatic field 
$\nabla\phi$. In particular, both entropy as well as internal energy separates 
electro-chemical-mechanical variables and the heat variable. We have 
\begin{align}
\sigma=\pl_e^{}\psi(e,c),\ \ \ \
\eta=-\pl_\theta^{}\psi(\theta),\ \ \ \ \mu=\pl_c^{}\psi(e,c),
\ \ \ \ \vecd=-\pl_{\vece}^{}\psi(\vece)
\end{align}
with $\eta$ the entropy and $\vecd$ the electric induction being equal to 
$\eps\vece$. The evolution variant of the most general system 
\eq{poro-steady-el-therm} then reads as: 
\begin{subequations}\label{poro-evol-temp+c+++}
\begin{align}\label{poro-evol-temp1+c+++}
&\hspace*{-2em}\varrho\DDT u-{\rm div}\,\sigma=f-z_{_{\rm DOP}}\nabla\phi&&\text{with }\ \ \ \sigma=\pl_e\psi_{_{\rm ME}}(e(u),c)
,&&
\\\label{poro-evol-temp2+c+++}
&\hspace*{-2em}\DT c-{\rm div}(\bbM(c,\theta)\nabla\mu)=r(c,\theta)&&\text{with }\ \ \ 
\mu=\pl_c^{}\psi_{_{\rm ME}}(e,c)+z\phi
,&&
\\\label{poro-evol-temp3+c+++}
&\hspace*{-2em}c_{\rm v}(\theta)\DT\theta-{\rm div}\big(\bbK(c,\theta)\nabla\theta\big)=
\bbM(c,\theta)\nabla\mu\Colon\nabla\mu+h(c,\theta)
-\mu\Cdot r(c,\theta),\hspace*{-12em}
\\\label{poro-evol-temp4+c+++}
&\hspace*{-2em}{\rm div}\big(\eps\nabla\phi)+z\Cdot c=z_{_{\rm DOP}}+{\rm div}(z_{_{\rm DOP}}u),
\hspace*{-3em}&&&&
\end{align}\end{subequations}
with the heat capacity $c_{\rm v}(\theta)=-\theta\psi_{_{\rm TH}}''(\theta)$, and with
the boundary conditions \eq{poro-steady-BC-thermo}. The heat production and 
chemical-reaction rates can be naturally assumed to vanish at $\theta=0$, 
i.e.\ $h(c,0)=0$ and $r(c,0)=0$, and then one can prove $\theta\ge0$.
Note that, instead of a full Maxwell electromagnetic system, we kept the 
electrostatic equation \eq{poro-evol-temp4+c+++}, which reflects the
well acceptable modelling assumption that the thermo-chemical-mechanical
processes are much slower than the electromagnetic processes and that
the electric currents are rather small so that the magnetic effects 
can be neglected. 
The balance of the electro-chemical-mechanical energy can be revealed by 
testing the particular equations (\ref{poro-evol-temp+c+++}a,b,d) successively
by $\DT u$, $\mu$, and $\DT\phi$:
\begin{align}\nonumber
&\hspace*{-4em}\frac{\d}{\d t}\bigg(\int_\Omega\frac\varrho2|\DT u|^2+
\psi_{_{\rm ME}}\big(e(u),c\big)
+z\phi
\Cdot c
-z_{_{\rm DOP}}\phi+z_{_{\rm DOP}}\nabla\phi\Cdot u\,\d x-\int_{\R^d}\frac\eps2|\nabla\phi|^2\,\d x\bigg)
\\[-.7em]&\hspace*{-4em}\nonumber\qquad\qquad
\qquad
+\!\!\lineunder{\int_\Omega
\bbM(c,\theta)\nabla\mu\Colon\nabla\mu\,\d x
+\int_{\Gamma}\alpha\mu^2\,\d S}{dissipation rate}\!
\\[-.2em]&\hspace*{-4em}\qquad\qquad\qquad
=\int_\Omega f\Cdot \DT u+r(c,\theta)\Cdot\mu\,\d x
+\int_{\GNeu}\!\!g\Cdot\DT u\,\d S
+\int_\Gamma\alpha\mu_{\rm ext}\mu\,\d S.
\label{poro-therm-engr-}
\end{align}
Using again 
\eq{calculus+}
and adding also \eq{poro-evol-temp3+c+++} tested by 1 we obtain the 
{\it total-energy balance}
\begin{align}\nonumber
&\hspace*{-4em}\frac{\d}{\d t}\bigg(\int_\Omega\!\!\!\!\!\linesunder{\frac\varrho2|\DT u|^2}{kinetic}{energy}\!\!\!\!\!\!
+
\!\!\!\!\!\linesunder{\psi_{_{\rm ME}}(e(u),c)+C_{\rm v}(\theta)_{_{_{_{_{}}}}}}{internal 
chemo-thermo-}{-mechanical energy}\!\!\!\!\!
\,\d x
+\int_{\R^d}\!\!\!\!\!\linesunder{\frac\eps2|\nabla\phi|^2}{electrostatic}{energy}\!\!\!\!\!\!\d x\bigg)
\\&\hspace*{-4em}\ =\!\!\!\!\!\linesunder{\int_{\GNeu}\!\!g\Cdot\DT u\,\d S+\int_\Omega\!f\Cdot \DT u+h(c,\theta)
\,\d x
}{power of mechanical load and of}{the heat from chemical reactions}\!\!\!\!\!
+\!\!\!\!\!\!\linesunder{\int_{\Gamma_{}}\!\alpha\mu_{\rm ext}\mu\,\d S}{power of}{chemical load}
\!\!\!\!\!+\!\!\!\!\!\linesunder{\int_{\Gamma_{}}\!\gamma(\theta_{\rm ext}{-}\theta)\,\d S}{power of the}{boundary heat flux}\!\!\!\!.
\label{poro-therm-engr}
\end{align} 
The heat part of the internal energy is 
$C_{\rm v}(\theta)=\psi_{_{\rm TH}}(\theta)-\theta\psi_{_{\rm TH}}'(\theta)$.
The analysis of the system \eq{poro-evol-temp+c+++} in general however does not
seem clear without modification by some gradient terms, as already mentioned at the end of Section~\ref{sec-FE}.
\end{remark}

\begin{remark}[{\sl Nonconvex free energies $\varphi$}]\upshape
Some applications call for nonconvex $\varphi$. E.g.\ a double-well potential
in $c$ is used in modelling of phase separation. The existence of the 
simplest static problem \eq{swelling-small-static-funct-u-c} obviously
requires some regularization, most conventionally by adding a gradient term like
$\frac12\epsilon|\nabla c|^2$. 
Yet, the nonconvexity certainly destroys the fixed-point argument in 
the ``non-static'' Propositions~\ref{prop-steady}--\ref{prop-steady-el-therm}.
The same comment applies for nonconvexity in terms of $e$ and, in particular,
arising from geometrical nonlinearity in a large-strain generalization like
in \cite{Anan11TMCT}. 
The possible nonexistence of steady states may be relatively 
expectable due to autonomous oscillations sometimes reported
in literature in the corresponding evolutionary systems, see e.g.\ 
\cite{BCMK05DRPE,SDHW10OPEM}, although strictly speaking 
those experimental experience rather documents nonexistence of stable
steady states only.
\end{remark}

\bigskip

{\small
\noindent{\it Acknowledgement}:
The author is very indebted for inspiring discussions with 
dr.\ J\"urgen Fuhrmann and prof.\ Alexander Mielke at Weierstrass Inst.\ 
Berlin about polymer-electrolyte fuel cells and about electroneutrality. 
}

\end{sloppypar}

\begin{thebibliography}{10}

\bibitem{Anan11TMCT}
L.~Anand.
\newblock A thermo-mechanically-coupled theory accounting for hydrogen
  diffusion and large elastic-viscoplastic deformations of metals.
\newblock {\em Intl. J. Solids Structures}, 48:962--971, 2011.

\vspace*{-.5em}\bibitem{BCMK05DRPE}
J.~Benziger, E.~Chia, J.~F. Moxley, and I.~G. Kevrekidis.
\newblock The dynamic response of {PEM} fuel cells to changes in load.
\newblock {\em Chemical Engineering Science}, 60:1743--1759, 2005.

\vspace*{-.5em}\bibitem{Biot41GTTS}
M.~A. Biot.
\newblock General theory of three-dimensional consolidation.
\newblock {\em J. Appl. Phys.}, 12:155--164, 1941.

\vspace*{-.5em}\bibitem{Fan52FPMT}
K.~Fan.
\newblock Fixed-point and minimax theorems in locally convex linear spaces.
\newblock {\em Proc. Nat. Acad. Sci. USA}, 38:121–126, 1952.

\vspace*{-.5em}\bibitem{Fuhr??MNMF}
J.~Fuhrmann.
\newblock Mathematical and numerical modeling of flow, transport and reactions
  in porous structures of electrochemical devices.
\newblock In P.~{Bastian et al.}, editor, {\em Simulation of Flow in Porous
  Media: Appl. in Energy and Environment}. J.Wiley, to appear.

\vspace*{-.5em}\bibitem{Latroche2004}
M.~Latroche.
\newblock Structural and thermodynamic properties of metallic hydrides used for
  energy storage.
\newblock {\em J. Physics \& Chemistry Solids}, 65:517--522, 2004.

\vspace*{-.5em}\bibitem{Libowitz1994}
G.G. Libowitz.
\newblock Metallic hydrides; fundamental properties and applications.
\newblock {\em J. Physics \& Chemistry Solids}, 55:1461--1470, 1994.

\vspace*{-.5em}\bibitem{ProWet09PEMF}
K.~Promislow and B.~Wetton.
\newblock {P}{E}{M} fuel cells: a mathematical overview.
\newblock {\em SIAM J. Appl. Math.}, 70:369--409, 2009.

\vspace*{-.5em}\bibitem{Raja07HAMF}
K.~R. Rajagopal.
\newblock On a hierarchy of approximate models for flows of incompressible
  fluids through porous solids.
\newblock {\em Math. Models Meth. Appl. Sci.}, 17:215--252, 2007.

\vspace*{-.5em}\bibitem{Roub13NPDE}
T.~Roub{\'\i}{\v{c}}ek.
\newblock {\em Nonlinear Partial Differential Equations with Applications}.
\newblock Birkh\"auser, Basel, 2nd edition, 2013.

\vspace*{-.5em}\bibitem{RouTom14THSM}
T.~Roub\'i\v{c}ek and G.~Tomassetti.
\newblock Thermomechanics of hydrogen storage in metallic hydrides: modeling
  and analysis.
\newblock {\em Discr. Cont. Dyn. Syst. B}, 14:2313--2333, 2014.

\vspace*{-.5em}\bibitem{SDHW10OPEM}
D.~G. Sanchez, D.~G. Diaz, R.~Hiesgen, I.~Wehl, and K.~A. Friedrich.
\newblock Oscillations of {PEM} fuel cells at low cathode humidification.
\newblock {\em J. Electroanalytical Chemistry}, 649:219--231, 2010.

\vspace*{-.5em}\bibitem{ShoSte04DPEM}
R.~E. Showalter and U.~Stefanelli.
\newblock Diffusion in poro-elastic media.
\newblock {\em Math. Methods Appl. Sci.}, 27:2131--2151, 2004.

\vspace*{-.5em}\bibitem{Stam65PDEE}
G.~Stampacchia.
\newblock Le probl\`eme de {D}irichlet pour les \'equations elliptiques du
  second ordre \`a coefficients discontinus.
\newblock {\em Ann. Inst. Fourier (Grenoble)}, 15:189--258, 1965.

\vspace*{-.5em}\bibitem{Visi84SCRR}
A.~Visintin.
\newblock Strong convergence results related to strict convexity.
\newblock {\em Comm. Partial Diff. Equations}, 9:439--466, 1984.

\vspace*{-.5em}\bibitem{Neum28TG}
J.~von Neumann.
\newblock Zur {T}heorie der {G}esellschaftsspiele.
\newblock {\em Math. Ann.}, 100:295--320, 1928.

\vspace*{-.5em}\bibitem{Zeid86NFAA}
E.~Zeidler.
\newblock {\em Nonlinear Functional Analysis and its Applications. {I}-{IV}}.
\newblock Springer, New York, 1985-90.

\end{thebibliography}
\end{document}